\documentclass[12pt]{amsart}
\usepackage{amsmath,amscd,amssymb,amsfonts}
\newcommand{\h}{\hbox}
\newcommand{\q}{\quad}
\newcommand{\nin}{\noindent}
\newcommand{\bs}{\par\bigskip}
\newcommand{\ms}{\par\medskip}
\newcommand{\sk}{\par\smallskip}
\newcommand{\bsn}{\par\bigskip\noindent}
\newcommand{\msn}{\par\medskip\noindent}

\newcommand{\ssb}{\raise.15ex\h{${\scriptscriptstyle\bullet}$}}
\newcommand{\ssc}{\,\raise.15ex\h{${\scriptstyle\circ}$}\,}
\newcommand{\msum}{\hbox{$\sum$}}
\newcommand{\mprod}{\hbox{$\prod$}}
\newcommand{\mcap}{\hbox{$\bigcap$}}
\newcommand{\mcup}{\hbox{$\bigcup$}}
\newcommand{\mopl}{\hbox{$\bigoplus$}}
\newcommand{\C}{{\mathbb C}}
\newcommand{\D}{{\mathbb D}}
\newcommand{\Q}{{\mathbb Q}}
\newcommand{\PP}{{\mathbb P}}
\newcommand{\R}{{\mathbb R}}
\newcommand{\RR}{{\mathbf R}}
\newcommand{\mm}{{\mathfrak m}}
\newcommand{\Z}{{\mathbb Z}}
\newcommand{\A}{{\mathcal A}}
\newcommand{\B}{{\mathcal B}}
\newcommand{\E}{{\mathcal E}}
\newcommand{\F}{{\mathcal F}}
\newcommand{\I}{{\mathcal I}}
\newcommand{\Hc}{{\mathcal H}}
\newcommand{\Lc}{{\mathcal L}}
\newcommand{\M}{{\mathcal M}}
\newcommand{\Oc}{{\mathcal O}}
\newcommand{\XX}{{\mathcal X}}
\newcommand{\Y}{{\mathcal Y}}
\newcommand{\Dt}{\widetilde{D}}
\newcommand{\jj}{\widetilde{j}}
\newcommand{\x}{\widetilde{x}}
\newcommand{\X}{\widetilde{X}}
\newcommand{\Yh}{\widehat{Y}}
\newcommand{\Xo}{{}\,\overline{\!X}{}}
\newcommand{\Gr}{{\rm Gr}}
\newcommand{\bl}{\bigl}
\newcommand{\br}{\bigr}
\newcommand{\into}{\hookrightarrow}
\newcommand{\simto}{\buildrel\sim\over\longrightarrow}
\newcommand{\onto}{\mathop{\rlap{$\to$}\hskip2pt\h{$\to$}}}
\newcommand{\muu}{\rlap{$\mu$}\hskip1.1pt\h{$\mu$}}
\newcommand{\indlim}{\rlap{\raise-5.5pt\h{$\,\to$}}{\rm lim}}

\newcommand{\ges}{\geqslant}
\newcommand{\les}{\leqslant}
\newcommand{\1}{\hskip1pt}
\begin{document}
\title[Weight zero part of first cohomology]
{Weight zero part of the first cohomology\\of complex algebraic varieties}
\author[M.~Saito]{Morihiko Saito}
\address{M. Saito : RIMS Kyoto University, Kyoto 606-8502 Japan}
\begin{abstract} We show that the weight 0 part of the first cohomology of a complex algebraic variety $X$ is a topological invariant, and give an explicit description of its dimension using a topological construction of the normalization of $X$, where $X$ can be reducible, but must be equidimensional. The first assertion is known in the $X$ compact case by A. Weber, where intersection cohomology is used. Note that the weight 1 or 2 part of the first cohomology is not a topological (or even analytic) invariant in the non-compact case by Serre's example.\end{abstract}
\maketitle
\centerline{\bf Introduction}
\bsn
Let $X$ be a complex algebraic variety which is assumed equidimensional and reduced in this paper. Let $\pi:\X\to X$ be the normalization. It is known that the underlying topological space of $\X$ can be constructed canonically from the underlying topological space of $X$, see \cite{GM1,McC} (and (2.2) below). Set
$$\F_X:={\rm Coker}(\pi^*:\Q_X\to\pi_*\Q_{\X}).$$
Then $\F_X$ is a constructible sheaf supported on the {\it non-unibranch} locus of $X$, and can be rather complicated in general. There is a long exact sequence of mixed Hodge structures
$$0\to H^0(X,\Q)\to H^0(\X,\Q)\to H^0(X,\F_X)\to H^1(X,\Q)\to H^1(\X,\Q)\to,
\leqno(1)$$
using the mapping cone construction in \cite{De3} or \cite{mhm} (where both give essentially the same, see \cite{mhc}). The following is known.
\msn
{\bf Proposition~1.} {\it The first cohomology $H^1(\X,\Q)$ has weight $\ges 1$.}
\ms
Indeed, this assertion follows from the surjectivity of the morphism of fundamental groups $\sigma_*:\pi_1(\XX)\to\pi_1(\X)$ for a resolution of singularities $\sigma:\XX\to\X$, using \cite[0.7 (B)]{FL} or \cite[Proposition 2.10]{Ko} (see \cite[Theorem~2.1]{ADH}), since it induces the injection of mixed Hodge structures $\sigma^*:H^1(\X,\Q)\into H^1(\XX,\Q)$.
\sk
We can also deduce Proposition~1 from the following (see (1.1) below).
\msn
{\bf Proposition~2.} {\it We have the injectivity of the canonical morphism of mixed Hodge structures}
$$H^1(\X,\Q)\into{\rm IH}^1(\X,\Q)={\rm IH}^1(X,\Q).
\leqno(2)$$
\ms
Indeed, it is well known that Proposition~1 still holds by replacing $H^1(\X,\Q)$ with the intersection cohomology ${\rm IH}^1(\X,\Q)$, see \cite[4.5.2]{mhm}. In this paper we prove the following.
\msn
{\bf Theorem~1.} {\it There are canonical isomorphisms}
$$\aligned W_0H^1(X,\Q)&={\rm Ker}\bl(H^1(X,\Q)\to{\rm IH}^1(X,\Q)\br)\\&={\rm Ker}\bl(H^1(X,\Q)\to H^1(\X,\Q)\br)\\&={\rm Coker}\bl(H^0(\X,\Q)\to H^0(X,\F_X)\br).\endaligned
\leqno(3)$$
\sk
(The first isomorphism is shown in the $X$ compact case by A.~Weber \cite[Theorem~1.7]{We}.)
\msn
{\bf Theorem~2.} {\it We have the equality}
$$\dim W_0H^1(X,\Q)=\dim H^0(X,\F_X)-b_0(\X)+b_0(X).
\leqno(4)$$
\ms
These assertions follow from the exact sequence (1) together with Propositions~1 and 2 if we have further the following.
\msn
{\bf Proposition~3.} {\it The cohomology group $H^0(X,\F_X)$ has weight $0$.}
\ms
This can be shown using mixed Hodge modules, for instance, see (1.2) below. Theorem~2 is well known in the $X$ curve case, where we have
$$\dim \F_{X,x}=r_{X,x}-1,$$
with $r_{X,x}$ the number of local irreducible components of $X$ at $x$. In the case $\dim X\ges 2$, it is not necessarily easy to calculate $\dim H^0(X,\F_X)$ explicitly, see (2.4) below. We can calculate it if $X$ satisfy certain strong conditions. Note that $\dim H^0(X,\F_X)$ can increase strictly by shrinking $X$, see Example~(2.5) below.
\sk
This work is partially supported by Kakenhi 15K04816. I thank N.~Budur for a good question about the topological invariance of $W_0H^1(X,\Q)$. This paper is written in order to answer his question.
\sk
In Section~1 we prove Theorems~1 and 2 by showing Propositions~2 and 3.
In Section~2 we review a topological construction of the normalization $\X$ of a complex analytic space $X$, and calculate $\dim H^0(X,\F_X)$ for some examples. In Appendix we review totally ramified cyclic coverings which are useful to construct examples of algebraic varieties with multi-branch points, and calculate the local monodromies of local branches..
\bs\bs
\vbox{\centerline{\bf 1. Proof of the main theorems}
\bsn
In this section we prove Theorems~1 and 2 by showing Propositions~2 and 3.}
\msn
{\bf 1.1.~Proof of Proposition~2.} In the notation of the introduction, let $K_{\X}^{\ssb}$ be the mapping cone of the canonical morphism $\Q_{\X}\to{\rm IC}_{\X}\Q[-n]$ so that we have the distinguished triangle
$$\Q_{\X}\to{\rm IC}_{\X}\Q[-n]\to K_{\X}^{\ssb}\buildrel{+1}\over\to\q\q\h{in}\q\q D^b_c(X,\Q),
\leqno(1.1.1)$$
where $n=\dim X$, and ${\rm IC}_{\X}\Q$ is the intersection complex, see \cite{BBD,GM2}.
\sk
Since $\X$ is normal, we have the following isomorphisms by an inductive construction of intersection complexes using iterated open direct images and truncations (see \cite[3.1]{GM2}):
$$\Hc^0\bl({\rm IC}_{\X}\Q[-n]\br)=\Q_{\X},\q\Hc^k\bl({\rm IC}_{\X}\Q[-n]\br)=0\,\,\,\,(\forall\,k<0).
\leqno(1.1.2)$$
By the long exact sequence of cohomology sheaves associated with (1.1.1), this implies
$$\Hc^kK_{\X}^{\ssb}=0\,\,\,\,(\forall\,k\les 0),\q\h{hence}\q H^0(\X,K_{\X}^{\ssb})=0,
\leqno(1.1.3)$$
since $\Gamma(X,*)$ is a left exact functor. So the injectivity of the first morphism in (2) follows from the long exact sequence of cohomology groups associated with (1.1.1):
$$H^0(\X,K_{\X}^{\ssb})\to H^1(\X,\Q)\to{\rm IH}^1(\X,\Q),
\leqno(1.1.4)$$
since ${\rm IH}^{\ssb}(\X,\Q)=H^{\ssb}(\X,{\rm IC}_{\X}\Q[-n])$ by definition.
\sk
The last isomorphism in (2) follows from the commutativity of the direct image by a {\it finite} morphism (in the strong sense) with the intermediate direct image (see Remark~(ii) below) in the cartesian diagram case, using the base change property as in \cite[4.4.3]{mhm}. Note that the direct image by a finite morphism $\pi$ is an exact functor of mixed Hodge modules (by the vanishing of $H^k\pi_*$ for $k\ne 0$). This finishes the proof of Proposition~2.
\msn
{\bf Remarks.} (i) In this paper, $a_X^*\Q\in D^b{\rm MHM}(X)$ with $a_X:X\to pt$ the structure morphism is denoted by $\Q_{h,X}$ for any complex algebraic variety $X$. Here $\Q$ denotes the trivial Hodge structure of type $(0,0)$ (see \cite{De1}), and ${\rm MHM}(X)$ is the category of mixed Hodge modules on $X$, see \cite{mhm}. We have $\Q_{h,X}[n]\in{\rm MHM}(X)$ if $X$ is smooth and purely $n$-dimensional.
\ms
(ii) The {\it intermediate direct image} of a mixed Hodge module $\M$ by an open embedding $j:U\into X$ is defined by
$$j_{!*}\M:={\rm Im}(H^0j_!\M\to H^0j_*\M)\in{\rm MHM}(X),$$
(see also \cite{BBD}). The intersection complex ${\rm IC}_X\Q_h$ in the category of mixed Hodge modules can be defined by $j_{!*}(\Q_{h,U}[n])$ for any smooth dense open subvariety $U\subset X$ assuming that $X$ is purely $n$-dimensional.
\msn
{\bf 1.2.~Proof of Proposition~3.} The constructible sheaf $\F_X$ underlies a bounded complex of mixed Hodge modules $\M^{\ssb}\in D^b{\rm MHM}(X)$ defined by
$$\M^{\ssb}:=C(\Q_{h,X}\to\pi_*\Q_{h,\X}).$$
Since $\pi$ is a finite morphism (in the strong sense), we have $H^ki_x^*\M^{\ssb}=0$ ($k\ne 0$) and $H^0i_x^*\M^{\ssb}$ is a Hodge structure of type $(0,0)$ for any $x\in X$, using the base change property as in \cite[4.4.3]{mhm}, where $i_x:\{x\}\into X$ is the inclusion.
\sk
For a sufficiently large number of points $x_k\in X$, we have the following injective morphism of mixed Hodge structures using the functorial morphisms ${\rm id}\to(i_x)_*i_x^*$ (see \cite[4.4.1]{mhm}):
$$H^0(X,\M^{\ssb}):=H^0(a_X)_*\M^{\ssb}\into\mopl_k\,H^0i_{x_k}^*\M^{\ssb}.
\leqno(1.2.1)$$
Here the injectivity is shown by using the finiteness of the cohomology $H^0(X,\F_X)$ underlying $H^0(X,\M^{\ssb})$. So Proposition~3 follows. (It seems also possible to apply the mapping cone construction in \cite{De3} to the base change of $\pi:\X\to X$ by $i_{x_k}:\{x_k\}\into X$. The details are left to the reader.)
\msn
{\bf 1.3.~Proof of Theorems~1 and 2.} Restricting the exact sequence (1) to the weight 0 part and using Propositions~1 and 3, we get the exact sequence
$$0\to H^0(X,\Q)\to H^0(\X,\Q)\to H^0(X,\F_X)\to W_0H^1(X,\Q)\to 0.
\leqno(1.3.1)$$
This implies Theorem~2 and the last equality in Theorem~1. The other equalities then follow from the exact sequence (1) and Proposition~2. This finishes the proof of Theorems~1 and 2.
\msn
{\bf 1.4.~Dependence of the weight filtration on the algebraic structure.} The weight filtration $W$ on $H^1(X,\Q)$ for a smooth {\it non-compact\1} surface is not even an {\it analytic\1} invariant. Indeed, $H^1(X,\Q)$ can be pure of weight either $1$ or $2$ depending on the algebraic structure in the case of Serre's example, where a 2-dimensional complex manifold $X$ has two smooth compactifications $\Xo_1$, $\Xo_2$ such that $\Xo_1$ is a ruled surface over an elliptic curve $E$ with divisor at infinity a section of the ruled surface, and $\Xo_2\cong\PP^2$ with divisor at infinity general three lines on $\PP^2$. More precisely we have in the first case
$$H^1(X_1,\Q)=\Gr^W_1H^1(X_1,\Q)=H^1(E,\Q),$$
where $X_i$ is the algebraic surface over $X$ determined by the compactification $\Xo_i$ ($i=1,2$) using GAGA, see \cite[Example 4.19]{PS}, \cite[Remark 2.12]{SS}, \cite[p.~232]{Ha1}, \cite[Appendix B, Example 2.0.1]{Ha2}.
\msn
{\bf Remark.} Using mixed 1-motives \cite[1.10]{De3}, it is rather easy to construct examples of smooth varieties $X$ such that $H^1(X,\Q)$ has weights $1$ and $2$, and the weight filtration $W_1H^1(X,\Q)$ as a submodule of $H^1(X,\Q)$ depends on the algebraic structure of $X$. For instance, let
$$H_{\Z}=\Z^3,\q H_{\C}=\C^3,\q F^0H_{\C}=\C\1v\subset H_{\C}\q\h{with}\q v:=(1,i,\alpha\1i)\in\C^3,$$
where $F^{-1}H_{\C}=H_{\C}$, $F^1H_{\C}=0$, and $\alpha$ is any irrational real number.
\sk
Let $\Gamma_{0,\Q}$ be any 1-dimensional $\Q$-vector subspace of $H_{\Q}:=H_{\Z}\otimes_{\Z}\Q$. Set
$$\Gamma:=H_{\Z},\q\Gamma_0:=\Gamma_{0,\Q}\cap\Gamma,\q\Gamma_1:=\Gamma/\Gamma_0.$$
Then $\Gamma_1$ is torsion-free, and $\Gamma_0\cap F^0H_{\C}=0$. Put
$$W_{-3}H_{\Z}=0,\q W_{-2}H_{\Z}=\Gamma_0,\q W_{-1}H_{\Z}=H_{\Z},$$
\vskip-2mm\nin
and
$$\aligned&V:=H_{\C}/F^0H_{\C},\q V_0:=\C\Gamma_0\subset V,\q V_1:=V/V_0,\\ &G:=V_0/\Gamma_0,\q X:=V/\Gamma,\q E:=V_1/\Gamma_1.\endaligned$$
Then $E$ is an elliptic curve. (Indeed, $(b+\alpha\1i)/(a+i)\notin\R$ for any $a,b,\in\Q$, since $\alpha$ is an irrational real number.) We have $G=\C^*$, and there is a fibration $X\to E$ whose fibers are parallel translates of $G\subset X$. This gives an {\it algebraic structure} as an semi-abelian variety on the complex Lie group $X$ (by compactifying the fibers). It corresponds to the mixed $\Z$-Hodge structure $H=\bl((H_{\Z},W),(H_{\C},F)\br)$ of weights $\{-2,-1\}$, which is the dual of $H^1(X,\Z)$, see {\it loc.\,cit}.
\bs\bs
\vbox{\centerline{\bf 2. Topological construction of the normalization.}
\bsn
In this section we review a topological construction of the normalization $\X$ of a complex analytic space $X$, and calculate $\dim H^0(X,\F_X)$ for some examples.}
\msn
{\bf 2.1.~Algebraic and analytic normalizations.} If a reduced complex algebraic variety $X$ is normal, then so is the associated analytic space $X^{\rm an}$; in particular, $X^{\rm an}$ is unibranch at any point, see for instance \cite{Mu}. This seems to be a rather nontrivial assertion (related to a version of Zariski's Main Theorem).
\msn
{\bf Remark.} The above theorem can be verified by using GAGA together with a well-known assertion that a reduced algebraic variety $X$ is normal if and only if the canonical morphism
$$\Oc_X\to\sigma_*\Oc_{\XX}$$
is an isomorphism for a resolution of singularities $\sigma:\XX\to X$ (and similarly for analytic spaces). Note that the isomorphism $\Oc_{X^{\rm an}}\simto\sigma_*\Oc_{\XX^{\rm an}}$ implies first that $X$ is unibranch at any point.
\msn
{\bf 2.2.~Local cohomology.} Let $X$ be a reduced complex analytic space of dimension $\les n$, where $n\ges 1$. It is well known that the the number of local irreducible components of dimension $n$ at $x\in X$ is given by
$$r_{X,x}:=\dim H^{2n}_{\{x\}}\Q_X,
\leqno(2.2.1)$$
which is defined purely {\it topologically} using the local cohomology.
\msn
{\bf Remarks.} (i) It is known that we have the isomorphism
$$H^{2n}_{\{x\}}\Q_X=H^{2n-1}(S_x\cap X,\Q),$$
where $S_x$ is a sufficiently small sphere with center $x$ in the ambient space containing $X$ locally, and is defined by using a distance function with $x$. This can be shown by using a topological cone theorem in the local singularity theory (where the theory of tubular neighborhood system \cite{Ma} is essential). This is usually used to show (2.2.1).
\ms
(ii) The assertion (2.2.1) can be proved by using only a resolution of singularities. Indeed, we have
$$H^{2n}_{\{x\}}\Q_X=H^{2n}_{\{x\}}j_!\Q_U,$$
where $j:U\into X$ is an open embedding such that $U$ is smooth and $Z:=X\setminus U\subset X$ is a closed analytic subset of dimension $<n$. This isomorphism follows by using the long exact sequence associated with the distinguished triangle
$$j_!j^{-1}\to{\rm id}\to i_*i^{-1}\buildrel{+1}\over\to,$$
where $i:Z\into X$ is the inclusion.
\sk
The assertion (2.2.1) is then reduced to the unibranch case, and can be proved by taking an embedded resolution of singularities
$$\sigma:(\XX,Y,E)\to(X,Z,x),$$
as follows. Here we may assume that $E:=\sigma^{-1}(x)$ is a divisor on $\XX$ (by taking first the blow-up along $x$), $Y:=\sigma^{-1}(Z)$ is a divisor with normal crossings on $\XX$, and $\sigma$ is an isomorphism over $U$. (Note that $E$ is also a divisor with normal crossings on $\XX$, since $E\subset Y$.) We have
$$\D\Q_U\cong\Q_U[2n],\q j_!\ssc\D=\D\ssc\RR\1j_*,\q i_x^!\ssc\D=\D\ssc i_x^*,\q (i_x)_*\ssc i_x^!=\RR\Gamma_{\{x\}},$$
with $i_x:\{x\}\into X$ the inclusion. Here $\D$ is the dual functor, see \cite{Ve}. The assertion (2.2.1) in the unibranch case is then reduced to
$$H^0i_x^*\RR\1j_*\Q_U=\Q.$$
Let $\jj:U\into\X$ be the inclusion so that $\sigma\ssc\jj=j$. By the base change property for the direct image by the proper morphism $\sigma:\X\to X$ and the pull-back by $i_x:\{x\}\into X$, the assertion is further reduced to
$$H^0(E,\RR\jj_*\Q_U|_E)=\Q,$$
and then follows from
$$\jj_*\Q_U|_E=\Q_E.$$
Here the {\it connectivity} of $E$ is essential.
\ms
(iii) The argument in Remark~(ii) above can be extended to the $\Z$-coefficient case. Here we have to show
$$H_{2n}(X,X\setminus\{x\};\Z)\,\bl(=H^0i_x^*\D\Z_X\br)=\Z,$$
in the unibranch case, see \cite[4.1]{GM1}. This is not necessarily equivalent to
$$H^{2n}(X,X\setminus\{x\};\Z)\,\bl(=H^{2n}_{\{x\}}\Z_X\br)=\Z,$$
since the $t$-structure on the derived category of bounded complexes of $\Z$-modules with finitely generated cohomologies is not self-dual (because of the problem of torsion), see \cite{BBD}.
\msn
{\bf 2.3.~Topological construction.} From now on, assume $X$ reduced and moreover purely $n$-dimensional (that is, $X$ is equidimensional with dimension $n$). The subset of $r$-branch points of $X$ can be defined by
$$X^{(r)}:=\{x\in X\mid r_{X,x}=r\}\subset X.$$
The local irreducible decomposition of $X$ is obtained by using the set of unibranch points $X^{(1)}$. (Indeed, for a sufficiently small open neighborhood $U_x$ of $x\in X$, it is enough to take the closures of the connected components of $X^{(1)}\cap U_x$.)
\sk
We will denote by $I_{X,x}$ the set of local irreducible components of $(X,x)$ for $x\in X$. (More precisely, $I_{X,x}$ is defined to be the inductive limit of the set $I(U_x,x)$ consisting of the connected components of $X^{(1)}\cap U_x$ whose closure contains $x$. Here $U_x$ runs over sufficiently small open neighborhoods of $x\in X$ such that $|I(U_x,x)|$ is independent of $U_x$.)
\sk
We can define a natural topology on
$$\X^{\rm top}:=\bigsqcup_{x\in X}I_{X,x}$$
so that an open neighborhood of $\x\in I_{X,x}$ is topologically identified with the normalization of the local irreducible component $\Gamma_{\x}$ of $X$ at $x$ corresponding to $\x\in I_{X,x}$ under the natural map
$$\pi^{\rm top}:\X^{\rm top}\to X.$$
Indeed, the local irreducible component $\Gamma_{\x}$ of $(X,x)$ defines a {\it multivalued} (set-theoretic) section of $\pi^{\rm top}$ over $\Gamma_{\x}\subset X$. (It is univalued only over the unibranch points of $\Gamma_{\x}$, for instance, if $X=\{y^2=x^2z\}\subset\C^3$.) We can get a fundamental neighborhood system of $\x\in\X^{\rm top}$ as the image of a fundamental neighborhood system of $(\Gamma_{\x},x)$ under this multivalued section (since the normalization is a {\it finite} morphism in the strong sense, that is, {\it proper and finite fibers}). This is a topological construction of the normalization $\X$ of a complex analytic space $X$, see \cite[4.1]{GM1}, \cite{McC} for {\it topological normalization}.
\msn
{\bf 2.4.~Calculation of $H^0(\X,\F_X)$.} With the notation of the introduction and (2.3), we have the isomorphism
$$\F_{X,x}={\rm Coker}\bl(\Q\into\mopl_{\x\in I_{X,x}}\,\1\Q\1\br)\q(x\in X),
\leqno(2.4.1)$$
and the monodromies of the local systems $\pi_*\Q_{\X}|_{S_i}$, $\F_X|_{S_i}$ are induced by the monodromy of the set $I_{X,x}$ ($x\in S_i$), where $S_i$ is a stratum of a Whitney stratification of $X$. In particular, the monodromy group of $\pi_*\Q_{\X}|_{S_i}$ is finite, and the monodromy action is semisimple.
\sk
Take $x_i\in S_i$ for each $i$. Let $\F^{\1\rm inv}_{X,x_i}\subset\F_{X,x_i}$ be the monodromy invariant part by the action of $\pi_1(S_i,x_i)$. Then
$$H^0(X,\F_X)=\mopl_i\,\F^{\1\rm inv}_{X,x_i}\,\,\,\h{if all the strata $S_i$ are closed subsets of $X$}.
\leqno(2.4.2)$$
\sk
In general we need further the information about the extensions between the local systems $\F_X|_{S_i}$. Choosing a path from $x_i$ to $x_j$ inside $S_j$ except for the end point, we get a morphism
$$\rho_{j,i}:\F_{X,x_i}\to\F_{X,x_j}\q\h{if}\q S_i\subset\overline{S}_j.$$
The composition $\rho_{k,j}\ssc\rho_{j,i}$ coincides with $\rho_{k,i}$ up to the monodromy action of $\pi_1(S_k,x_k)$ if $S_i\subset\overline{S}_j\subset\overline{S}_k$. We have
$$H^0(X,\F_X)=\bl\{(\xi_i)\in\mopl_i\,\F^{\1\rm inv}_{X,x_i}\mid\rho_{j,i}\1\xi_i=\xi_j\,\,\,\h{if}\,\,\,S_i\subset\overline{S}_j\br\}.
\leqno(2.4.3)$$
\msn
{\bf Remark.} It is not necessarily easy to construct an example of irreducible projective variety $X$ such that
$$\dim H^0(X,\F_X)>0\q\h{with}\q\dim X>1,$$
except for the case where $X=C\times Z$ with $\dim C=1$ or a quotient of $C\times Z$ by a finite group action.
For instance it is known that, if $X\subset\PP^N$ is a {\it global complete intersection} with $\dim X>1$, then
$$H^1(X,\Q)=0.$$
Here Artin's theorem in \cite{BBD} is used, see also \cite{Di}. (The latter theorem can be deduced also from Cartan's Theorem~B together with the Riemann-Hilbert correspondence, see for instance \cite[Lemma~2.1.18]{mhp}.
Note that the shifted constant sheaf $\Q_X[\dim X]$ (or its dual) for a {\it locally complete intersection} $X\subset Y$ corresponds to a regular holonomic $D$-module on a complex manifold $Y$ by the Riemann-Hilbert correspondence. This can be verified easily by using algebraic local cohomology, see for instance a remark after \cite[1.5.1]{RSW}.)
\msn
{\bf 2.5.~Examples.} (i) For integers $a,b,c\ges 2$, set
$$X:=\{y^b=x^b(x^a+z^c)\}\subset\C^3.$$
This is irreducible at 0. Indeed, we get the normalization
$$\X:=\{v^b=x^a+z^c\}\subset\C^3,$$
by the blowing up along the center $\{x=y=0\}\subset\C^3$, where $v=y/x$ (for the normality of hypersurfaces with isolated singularities, see for instance \cite[II, Proposition 8.23]{Ha2}).
\sk
The non-unibranch locus coincides with the $b$-branch locus
$$X^{(b)}=S:=\{x=y=0\}\setminus\{0\}\subset X.$$
This is a stratum of the Whitney stratification. Note that, restricting $\X$ to $\{x=0\}$, we get
$$\pi^{-1}(\overline{S})=\{v^b=z^c\}.$$
\sk
Since $S$ is not closed in $X$, and $\F_X$ is the 0-extension of $\F_X|_S$, we have
$$H^0(X,\F_X)=0.
\leqno(2.5.1)$$
In this example $H^1(X,\Q)$ also vanishes by the exact sequence (1), since $\X$ is contractible so that $H^1(\X,\Q)=0$.
\ms
(ii) In the above example, the action of the monodromy on $I_{X,s}$ by a generator of $\pi_1(S,s)$ ($s\in S$) is identified with
$$\Z/b\Z\ni k\mapsto k+c\in\Z/b\Z,$$
using the description of $\pi^{-1}(\overline{S})$ in (i). One may also use the description of local irreducible components of $X$ given by
$$y=x(x^a+z^c)^{1/b},$$
where $x^a$ can be neglected essentially (by restricting to $\{|x|^a\ll|z|^c\}\subset\C^3$), see also (A.4.1) below.
\sk
The number of orbits by the monodromy action on $I_{X,s}$ is then given by
$$e:={\rm GCD}(b,c),$$
since $\Z/(b\Z+c\Z)=\Z/e\Z$.
\sk
Set $X':=X\setminus\{0\}$ (or $X':=X\setminus C$ with $C\subset X$ any curve containing $0$ and not contained in $\overline{S}$). Since the monodromy action on $\pi_*\Q_{\X}|_S$ is semisimple, we then get by (2.4.1--2)
$$\dim H^0(X',\F_{X'})=e-1.
\leqno(2.5.2)$$
\msn
{\bf Remark.} In the case $X':=X\setminus\{0\}$, we can calculate $H^1(\X',\Q)$ as in theory of isolated hypersurface singularities \cite{Mi} (using the Wang sequence), and get
$$H^1(\X',\Q)\cong H^2(F_{\!f},\Q)_1(1).
\leqno(2.5.3)$$
Here $\X':=\X\setminus\{0\}$ with $\X=f^{-1}(0)\subset\C^3$ for $f:=v^b-x^a-z^c$, and the right-hand side is the unipotent monodromy part of the Milnor fiber cohomology of $f$ up to a Tate twist. (Note that the Milnor monodromy is semisimple so that $N=0$ in the weighted homogeneous isolated singularity case, where $N:=(2\pi i)^{-1}\log T_u$ with $T_u$ the unipotent part of the Milnor monodromy $T$.) Indeed, we have the canonical isomorphism
$$i_x^*\RR(j_U)_*\Q_U=C\bl(N:i_x^*\psi_{f,1}\Q_V\to i_x^*\psi_{f,1}\Q_V(-1)\br)[-1],$$
where $U:=\C^3\setminus f^{-1}(0)$ with the inclusion $j_U:U\into V:=\C^3$, and $\psi_{f,1}$ is the unipotent monodromy part of the nearby cycle functor $\psi_f$ \cite{De2}. This implies
$$H^4i_x^*i_{\X}^!\Q_V=H^3i_x^*\RR(j_U)_*\Q_V=H^2(F_{\!f},\Q)_1(-1),$$
where $i_{\X}:\X\into V=\C^3$ is the inclusion. Since $\D\Q_V=\Q_V(3)[6]$, we have
$$\D(H^4i_x^*i_{\X}^!\Q_V)=H^{-4}i_x^!i_{\X}^*\Q_V(3)[6]=H^2i_x^!\Q_{\X}(3)=H^1(\X',\Q)(3),$$
where the last isomorphism follows from the distinguished triangle
$$(i_x)_*i_x^!\to{\rm id}\to\RR(j_{\X'})_*j_{\X'}^*\buildrel{+1}\over\to,$$
with $j_{\X'}:\X'\into\X$ the inclusion. (Here $i_x$ denotes also the inclusion into $\X$.) So (2.5.3) follows from the above two isomorphisms together with the self-duality isomorphism
$$\D\bl(H^2(F_{\!f},\Q)_1\br)=H^2(F_{\!f},\Q)_1(3),$$
which is a consequence of the duality for the unipotent monodromy part of the vanishing cycle functor $\varphi_{f,1}$, see for instance \cite[5.2.3]{mhp}. (Here the intersection form cannot be used on this unipotent monodromy part, see \cite{St1,St2}.)
\sk
The unipotent monodromy part of the vanishing cohomology $H^2(F_{\!f},\Q)_1$ has pure Hodge structure of weight 3, since $N=0$, see \cite{St1,St2} and also \cite[5.1.6]{mhp}. Its dimension is given by the sum of the coefficients of $t$ and $t^2$ in the Steenbrink spectrum as in \cite{St2}, and the latter is expressed in this case by
$${\rm Sp}(f)=\biggl(\frac{t^{1/a}-t}{1-t^{1/a}}\biggr)\biggl(\frac{t^{1/b}-t}{1-t^{1/b}}\biggr)\biggl(\frac{t^{1/c}-t}{1-t^{1/c}}\biggr),$$
see for instance \cite[1.9]{wh} and the references there. (In the case $a=b=c$, we may also use the Thom-Gysin sequence (see for instance \cite[1.3]{RSW}) together with a well-known theorem of Griffiths on rational integrals for the complement of a smooth curves on $\PP^2$.) By the symmetry of the coefficients of the spectrum, we then get
$$\dim H^1(\X',\Q)=2\cdot\#\bl\{(i,j,k)\in\Z_{>0}^3\mid\tfrac{i}{a}+\tfrac{j}{b}+\tfrac{k}{c}=1\br\}.$$
This description implies that we may have
$$W_1H^1(X',\Q)=0,\q\h{but}\q H^1(X',\Q)=H^1(\X')\ne 0,$$
in the case $e={\rm GCD}(b,c)=1$ and $\tfrac{a}{b},\tfrac{a}{c}\in\Z$, for instance, if $(a,b,c)=(6,3,2)$. (Note that $H^1(X',\F_{X'})=0$, since $\pi^{-1}(X'{}^{(b)})\cong\C^*$ by $e=1$.)
\bs\bs
\centerline{\bf Appendix. Totally ramified cyclic coverings}
\bsn
In this Appendix we review totally ramified cyclic coverings which are useful to construct examples of algebraic varieties with multi-branch points, and calculate the local monodromies of local branches.
\msn
{\bf A.1.~Construction.} Let $Y$ be a complex algebraic variety, $D$ be a locally principal effective divisor on $Y$, and $\Lc$ be an invertible sheaf (that is, a locally free sheaf of rank 1) on $Y$ such that
$$\iota:\Lc^{\otimes\,m}\simto\Oc_Y(-D)\subset\Oc_Y,
\leqno{\rm(A.1.1)}$$
where $m\ges 2$. This isomorphism gives a structure of an $\Oc_Y$-algebra on
$$\A_Y:=\mopl_{i=0}^{m-1}\,\Lc^{\otimes\,i}=\B_Y/\I(D,\Lc)\q\h{with}\q\B_Y:=\mopl_{i\ges 0}\,\Lc^{\otimes\,i},$$
where $\I(D,\Lc)\subset\B_Y$ is an ideal sheaf locally generated over $\B_Y$ by
$$f_u:=\otimes^mu-\iota(\otimes^mu)\in\Lc^{\otimes\,m}\oplus\Oc_Y\subset\B_Y,$$
for a local generator $u\in\Lc$.
\sk
The {\it totally ramified cyclic covering\1} of $Y$ of degree $m$ associated with $(D,\Lc)$ is then defined by
$$X(Y,D,\Lc,m):=({\rm Spec}_Y\A_Y)_{\rm var}\,\subset\,L^{\vee}=({\rm Spec}_Y\B_Y)_{\rm var}.$$
Here $(S)_{\rm var}$ is the algebraic variety associated with a scheme $S$ of finite type over $\C$ in general, and is defined by the set of closed points of $S$ together with the induced topology and structure sheaf. We will denote $X(Y,D,\Lc,m)$ also by $X$ to simplify the notation.
\sk
This construction is quite well known in algebraic geometry and Hodge theory. It was, for instance, used to generalize Kodaira's vanishing theorem by applying C.P.~Ramanujam's idea \cite{Ram} so that the vanishing theorem is reduced to the weak Lefschetz type theorem using Hodge theory, see \cite[p.~151]{Na}, \cite[Proof of Proposition 2.33 and Remark 2.34(2)]{mhm} (and also \cite{Vi} where the {\it normalization} cannot be taken in our case since the multi-branch points disappear as in Example~(2.5)(i)). In the literature it does not seem that the {\it monodromies of local branches} are studied.
\sk
Since we assume that $\Lc$ is an invertible sheaf and $D$ is a locally principal divisor, the construction is quite trivial. By definition, $X$ is a locally principal divisor on the dual line bundle $L^{\vee}$ of the line bundle $L$ corresponding to the locally free sheaf $\Lc$ (that is, the sheaf of local sections of $L$ is $\Lc$). Choosing a local trivialization of $L$ (or equivalently, of $L^{\vee}$), $X$ is locally defined by
$$t^m=g_u\q\h{in}\q L^{\vee},
\leqno{\rm(A.1.2)}$$
where $g_u:=\iota(\otimes^mu)\in\Oc_Y$ for a local generator $u\in\Lc$ giving the local trivialization of $L^{\vee}$, and $t$ is the coordinate of the fiber of the line bundle $L^{\vee}$ associated with the local trivialization of $L^{\vee}$. (By definition $\Lc\subset\B_Y$ is identified with functions on $L^{\vee}$ which are {\it linear} on each fiber of $L^{\vee}$, and $t$ corresponds to $u$ by this identification.) Note that $g_u$ is a local defining function of $D$ by (A.1.1), and $X$ is a ramified finite cyclic Galois covering of $Y$, where a generator of the covering transformation group acts by multiplication by $\exp(2\pi\sqrt{-1}/m)$ on the variable $t$. We say that $X$ is {\it totally ramified,} since $\pi_Y^{-1}(D)_{\rm red}$ is isomorphic to $D_{\rm red}$, where $\pi_Y:X\to Y$ is the canonical morphism.
\bsn
{\bf Remarks.} (i) If we fix a divisor $D$, then the invertible sheaf $\Lc$ satisfying (A.1.1) is unique {\it up to} a tensor product with an invertible sheaf $\E$ satisfying $\E^{\otimes\,m}=\Oc_Y$ (which corresponds to an $m$-torsion point of the Picard group ${\rm Pic}(X)=H^1(X,\Oc_X^*)$). For a fixed divisor $D$, the totally ramified cyclic covering $X$ of degree $m$ associated with $(D,\Lc)$ depends (as an algebraic variety over $Y$) on the choice of the invertible sheaf $\Lc$ satisfying (A.1.1) in general. This can be seen, for instance, in the case of $m=2$ by looking at the cokernel of the canonical injection $\Oc_Y\into(\pi_Y)_*\Oc_X$, since the latter is isomorphic to $\Lc$ by the definition of $\Oc_X$ (that is, $(\pi_Y)_*\Oc_X=\A_Y$).
\ms
(ii) The construction of the totally ramified cyclic covering of degree $m$ associated with $(D,\Lc)$ is {\it compatible with base changes by dominant morphisms of $Y$}. More precisely, for a morphism $\phi:Y'\to Y$ such that the image of any irreducible component of $Y'$ is not contained in $D$, there is the canonical isomorphism
$$X(Y',\phi^*D,\phi^*\Lc,m)=X(Y,D,\Lc,m)\,{\times}_Y\,Y'.
\leqno{\rm(A.1.3)}$$
\msn
{\bf A.2.~Case $Y$ smooth.} In the notation of (A.1), assume $Y$ connected and {\it smooth\1} so that any divisor on $Y$ is locally principal. We have the irreducible decomposition
$$D=\msum_{i\in I}\,a_iD_i\q(a_i\in\Z_{>0}),$$
For $y\in D$, set
$$r_y(D,m):={\rm GCD}\bl(m,a_i\,(i\in I_y)\br)\q\h{with}\q I_y:=\{i\in I\mid y\in D_i\}.$$
In the notation of (2.3), we can prove the following isomorphism as topological spaces
$$\pi_Y:X^{(r)}\simto D(m,r):=\{y\in D\mid r_y(D,m)=r\},
\leqno{\rm(A.2.1)}$$
which is induced by the isomorphism as reduced algebraic varieties
$$\pi_Y:(X|_D)_{\rm red}\simto D_{\rm red},
\leqno{\rm(A.2.2)}$$
where the left-hand side means the reduced restriction of $X$ over $D_{\rm red}$. Indeed, we have a factorization of a defining function of $X$ in the line bundle
$$t^m-g=\mprod_{q=1}^e\,(t^{m/e}-\theta^{\,q}h)\q\h{if}\q g=h^e\q\h{with}\q m/e\in\Z,
\leqno{\rm(A.2.3)}$$
where $\theta=\exp(2\pi\sqrt{-1}/e)$. So the isomorphism (A.2.1) is reduced to the assertion that $X$ is unibranch at a point over $y\in D$ if $r_y(D,m)=1$.
\sk
The latter assertion is verified by taking locally a sufficiently general line $C$ in $Y$ which is sufficiently close to $y$ and intersecting only smooth points of the $D_i$ ($i\in I_y$), and then restricting the covering over $C$. Indeed, it is well known (and is easy to show) that the local monodromies of the unramified cyclic covering over the complement of $D$ are given around a intersection point of $C$ and $D_i$ by $\gamma^{-a_i}$. Here $X|_C$ is locally defined by $t^m=z^{a_i}$ in $\C\times C$ with $z$ a local coordinate of $C$, and we have the factorization
$$t^m-z^{a_i}=\mprod_{q=1}^m(t-\theta^{\prime\,q}\1z^{a_i/m})$$
over any simply connected open subset of $C\setminus D$ with $\theta':=\exp(2\pi\sqrt{-1}/m)$. Note that $\gamma$ is identified with the generator $\exp(2\pi\sqrt{-1}/m)$ of the covering transformation group $\muu_m$ of the cyclic covering $X\to Y$, where
$$\muu_m:=\{\theta\in\C^*\mid\theta^m=1\},
\leqno{\rm(A.2.4)}$$
and the latter acts by the natural multiplication on the variable $t$ in (A.1.2).
(Recall that the cyclic covering is an abelian covering so that the local monodromy is independent of a path from a base point to the intersection point inside $C$, since it is unique up to a conjugation.) So the restriction of the covering over $C\setminus D$ is connected if $r_y(D,m)=1$.
\msn
{\bf Remark.} If ${\rm GCD}\bl(m,a_i\,(i\in I)\br)=1$, then $X$ is irreducible by an argument similar to the above one by taking sufficiently general curves on $Y$, at least if $Y$ is quasi-projective. (Note that $Y$ is assumed connected.) Its converse does not necessarily hold as a consequence of Remark after (A.1) unless $H_1(Y,\Z)=0$.
\msn
{\bf A.3.~Calculation of local monodromies of local branches.} With the notation and assumption of (A.2), we have a stratification of $Y$ such that its strata (which are not assumed smooth) are given by
$$\aligned&Y_{J,p}^{\circ}:=Y_{J,p}\setminus\bl(\,\mcup_{i\notin J}\,D_i\br),\\ \h{with}\q\q&Y_{J,p}\,(k\in K_J)\q\h{irreducible components of}\,\,\,\,\mcap_{i\in J}\,D_i.\endaligned$$
These are indexed by $J\subset I$, $k\in K_J$ satisfying the condition:
$$\bl\{i\in I\mid D_i\supset Y_{J,p}\}=J.$$
We fix a non-empty $J\subset I$ and $k\in K_J$ satisfying this condition (where $Y_{J,p}\subset D_{\rm red}$). We have the local system of finite sets on $Y_{J,p}^{\circ}$
$$\{I_{X,\pi_Y^{-1}(y)}\}_{y\,\in\,Y_{J,p}^{\circ}},$$
consisting of local branches, where $Y_{J,p}^{\circ}$ is identified with $\pi_Y^{-1}Y_{J,p}^{\circ}$ by the isomorphism (A.2.2). Here a local system of finite sets means a locally constant sheaf with stalks finite sets. (We will see that it is a principal $\muu_e$-bundle at the end of this section.) We can describe its local monodromies around $Y_{J,p}\setminus Y_{J,p}^{\circ}$ as follows.
\sk
The stratum $Y_{J,p}^{\circ}$ and its closure $Y_{J,p}$ may have singularities in general. In the singular case, we apply Hironaka's resolution of singularities to the variety $Y_{J,p}$ inside $Y$. It is given by iterating blowing-ups along smooth centers contained in the singular locus of the proper transform of $Y_{J,p}$. We then get a projective morphism of smooth algebraic varieties
$$\sigma_{\Y}:\Y\to Y,$$
inducing an isomorphism over the complement of the singular locus of $Y_{J,p}$ in $Y$, and such that the proper transform $\Y_{J,p}$ of $Y_{J,p}$ in $\Y$ is smooth. (This procedure is unnecessary in the $Y_{J,p}$ smooth case.)
\sk
We further iterate blowing-ups along smooth centers contained in the non-normal-crossing locus of the intersection of the proper transform $\Y_{J,p}$ with the the total transforms of $\mcup_{i\notin J}\,D_i$ so that this intersection becomes a divisor with normal crossings on $\Y_{J,p}$. (It is not required that the total transforms of $\mcup_{i\notin J}\,D_i$ is a divisor with normal crossings on $\Y$.) This is necessary to calculate the local monodromies, since we can do it only at smooth points of the reduced intersection of $\Y_{J,p}$ and the total transforms of $\mcup_{i\notin J}\,D_i$. (Here reduced intersection means that we consider the associated reduced variety. This is needed since we talk about smooth points.)
\sk
Let
$$Z=\msum_k\,m_kZ_k$$
be the restriction of the total transforms of the divisor $\msum_{i\notin J}\,a_iD_i$ to $\Y_{J,p}$, where the $Z_k$ are the irreducible components of the divisor with normal crossings (by the above construction). Set
$$\XX:=X\times_Y\Y,\q\h{and}\q\Y_{J,p}^{\circ}:=\Y_{J,p}\setminus Z.$$
Note that $\XX$ is a totally ramified covering of $\Y$, see Remark~(ii) after (A.1). Using the isomorphism (A.2.2) for $\pi_{\Y}:\XX\to\Y$, the local monodromy of $\{I_{\XX,\pi_{\Y}^{-1}(y)}\}_{y\in\Y_{J,p}^{\circ}}$ around a general point $y\in Z_k\setminus\bl(\mcup_{k'\ne k}\,Z_{k'}\br)$ is given by
$$\gamma^{-m_k},
\leqno{\rm(A.3.1)}$$
where the generator $\gamma$ of the covering transformation group $\muu_m$ of $\pi_{\Y}:\XX\to\Y$ can be identified with a generator of the covering transformation group of $\widetilde{\XX}|_{\Y_{J,p}^{\circ}}\to\Y_{J,p}^{\circ}$ (with $\widetilde{\XX}$ the normalization of $\XX$). Indeed, the assertion can be reduced to the case $\dim\Y_{J,p}=1$ and $Z_k$ is a point, by restricting to a sufficiently general smooth subvariety of $\Y$ which is transversal to $Z_k$. Here we can further iterate blowing-ups along $Z_k$ so that the proper transform of $D_i$ for $i\notin J$ does not contain $Z_k$, and only the exceptional divisor of the last blow-up contains $Z_k$. (Here we may also assume that the proper transforms of local branches of $D_i$ for $i\in J$, which intersect $Y_{J,p}$ without containing it completely, do not contain $Z_k$.)
\sk
Set
$$e:={\rm GCD}\bl(m,a_i\,(i\in J)\br).$$
On a neighborhood of $Z_k$, we then get the factorization of a defining function of $\XX$ in the line bundle
$$t^m-z^{m_k}\,\mprod_{i\in J}\,h_i^{a_i}=\mprod_{q=1}^e\,\bl(t^{m/e}-\theta^{\,q}\,z^{m_k/e}\,\mprod_{i\in J}\,h_i^{a_i/e}\br).
\leqno{\rm(A.3.2)}$$
Here $\theta:=\exp(2\pi\sqrt{-1}/e)$, $h_i$ is the pull-back of a local defining function of $D_i$, and $z$ is an appropriate coordinate defining the exceptional divisor of the last blow-up along $Z_k$. (Note that $m/e$ and the $a_i/e$ ($i\in J$) are integers by the definition of $e$, but $m_k/e$ is not necessarily.) So we get the assertion (A.3.1), since $\gamma$ acts by multiplication by $1^{1/m}$ on the variable $t$, and hence by $1^{1/e}$ on $t^{m/e}$, where $1^{1/k}:=\exp(2\pi\sqrt{-1}/k)$ for $k\in\Z_{>0}$.
\sk
The above argument shows that the local branches form a principal $\muu_e$-bundle over $\Y_{J,p}^{\circ}$, which is identified with a locally constant sheaf with stalks finite sets having a transitive free action of $\muu_e$ (where $\muu_e$ is as in (A.2.4)).
\msn
{\bf Remark.} The global monodromies of local branches are determined by local monodromies if $H_1(\Y_{J,p},\Z)=0$. Indeed, for two principal $\muu_e$-bundles $P_1$, $P_2$ over $\Y_{J,p}^{\circ}$, $\Hc om_{\,\muu_e}(P_1,P_2)$ is also a principal $\muu_e$-bundle over $\Y_{J,p}^{\circ}$. The latter can be extended over $\Y_{J,p}$ if the local monodromies of $P_1$,$P_2$ coincide. Here the local monodromies are identified with morphisms
$$H_1(U_y\cap\Y_{J,p}^{\circ},\Z)\to\muu_e,$$
with $U_y$ an open neighborhood of $y$ in $\Y_{J,p}$ (since the target is abelian). A similar assertion holds for global monodromies with $U_y\cap\Y_{J,p}^{\circ}$ replaced by $\Y_{J,p}^{\circ}$.
\msn
{\bf A.4.~$2$-dimensional case.} With the notation and assumption of (A.3), assume further $\dim Y=2$, and $Y_{J,p}=D_i$. In this case $\Y$ is obtained by iterating point-center blowing-ups, and $Z_k$ is a point of $\Y_{J,p}$ (which is a proper transform of $D_i$), and is identified with a point $y'$ of the normalization $\Dt_i$ of $D_i$ corresponding to an analytic branch $(D'_i,y)\subset(D_i,y)$. Here $e={\rm GCD}(m,a_i)$. The multiplicity $m_k$ of $Z_k$ in $Z$ (see (A.3)) can be given by
$$m_{y'}:=\msum_{j\ne i} a_j(D'_i,D_j)_y,
\leqno{\rm(A.4.1)}$$
with $(D'_i,D_j)_y$ the intersection number of $D'_i$ and $D_j$ at $y$, and we have the equalities
$$(D'_i,D_j)_y:=\dim_{\C}\Oc_{Y^{\rm an},y}/(g'_i,g_j)=\#\bl(g_i^{\prime\,-1}(0)\cap g_j^{-1}(c)\br).
\leqno{\rm(A.4.2)}$$
Here $g'_i$, $g_j$ are (reduced) local defining holomorphic functions of $D'_i$, $D_j$ around $y$, and $c$ is a nonzero complex number with $|c|$ sufficiently small so that $g_i^{\prime\,-1}(0)$ and $g_j^{-1}(c)$ intersect transversally at smooth points. (Here we cannot use the finite determinacy of hypersurface isolated singularities, since we have to treat two functions $g'_i,g_j$ at the same time, although $g_j$ is algebraic in the initial coordinates.) As for the middle term of (A.4.2), we have the canonical isomorphism
$$\Oc_{Y^{\rm an},y}/(g'_i,g_j)=\Oc_{D'_i,y}\otimes_{\Oc_{Y^{\rm an},y}}\Oc_{D^{\rm an}_j,y}.$$
If $D_i$ is {\it unibranch} (where $D'_i=D_i$ and $g'_i=g_i$ are algebraic), the first equality of (A.4.2) is a theorem (in the proper intersection case), see \cite[Section 7.1]{Fu} (and also \cite{Se}).
\sk
The last term of (A.4.2) coincides with the multiplicity $m_{i,j}$ of $Z_k$ in the restriction of the total transform of $D_j$ to the proper transform of $D_i$, where the latter is identified with the normalization $\Dt_i$. Note that there is a local coordinate $\widetilde{z}$ of $(\Dt_i,y')$ such that the pull-back of $g_j$ to the proper transform $(\Dt_i,y')$ by the composition $\Dt_i\into \Y\to Y$ coincides with $\widetilde{z}^{\,m_{i,j}}$. (This can be used to show (A.4.2) in the $D_i$ unibranch case, see Remark~(ii) below.)
\sk
To show (A.4.2), consider the morphism
$$\rho:(\C^2,0)\ni(y_1,y_2)\mapsto(z_1,z_2)=(g'_i(y_1,y_2),g_j(y_1,y_2))\in(\C^2,0),$$
where $(y_1,y_2)$ is an analytic local coordinate system of $(Y,y)$. The assertion is then shown by proving the {\it finiteness} and {\it flatness} of the morphism
$$\rho|_U:U\to U'\q\h{with}\q U:=\rho^{-1}(U'),
\leqno{\rm(A.4.3)}$$
for a sufficiently small open neighborhood $U'$ of $0\in\C^2$ satisfying
$$g_i^{\prime\,-1}(0)\cap g_j^{-1}(0)\cap U=\{0\},\q\h{that is,}\q(\rho|_U)^{-1}(0)=\{0\}.$$
Indeed, finiteness follows from Weierstrass preparation theorem (using a graph embedding together with the factorization $\C^4\to\C^3\to\C^2$). So the direct image $\rho_*\Oc_U$ is a coherent $\Oc_{U'}$-module, and hence $R:=\C\{y_1,y_2\}$ is a finite $R'$-module with $R':=\C\{z_1,z_2\}$.
\sk
For flatness, consider a short exact sequence of finite $R'$-modules
$$0\to K\to P\to R\to 0,$$
where $P$ is free, $P/\mm'P\to R/\mm'R$ is an isomorphism with $\mm':=(z_1,z_2)\subset R'$ the maximal ideal, and $K$ is defined by the short exact sequence. (Note that the surjectivity of $P\to R$ follows from the isomorphism modulo $\mm'$ using Nakayama's lemma.) Since $g'_i,g_j$ form a regular sequence in $R=\C\{y_1,y_2\}$, the Koszul complex defined by them vanishes except for the highest degree, see \cite{Ei,Se}. This implies that
$${\rm Tor}^{R'}_1(R'/m'R',R)=0,$$
since the left-hand side coincides with the middle term of the Koszul complex for $g'_i,g_j$. Using the long exact sequence associated with the above short exact sequence, we then get
$$K/\mm'K=0,\q\h{hence}\q K=0\q\h{(by Nakayama's lemma).}$$
Thus $R$ is finite free over $R'$, see also \cite[Theorem 6.8]{Ei}. (Note that freeness is equivalent to flatness under the finiteness hypothesis by an argument similar to the above one, see \cite{Ei,Se}.) So the coherent sheaf $\rho_*\Oc_U$ is a finite free $\Oc_{U'}$-module (shrinking $U'$ if necessary). We then get the last equality of (A.4.2), since both sides coincide with the rank of the free $\Oc_{U'}$-module $\rho_*\Oc_{U}$ counted at $(0,0)$ and $(0,c)\in\C^2$ respectively.
\msn
{\bf Remarks.} (i) We can show the finiteness and flatness of (A.4.3) using algebraic geometry in the $D_i$ {\it unibranch} case, although Zariski topology is not fine enough for the assertion to hold with $U'$ a Zariski open neighborhood of $0\in\C^2$. By Grothendieck's version of Zariski's main theorem (see \cite[Theorem 4.4.3]{Gr2}), we can extend the germ of a morphism
$$(g_i,g_j):(Y,y)\to(\C^2,0)\q\h{(with $g_i,g_j$ algebraic)}$$
to a finite morphism $\Yh\to\C^2$, where $\Yh$ contains an open subvariety $Y^{\circ}$ of $Y$ as a dense open subvariety, and $\Yh\setminus Y^{\circ}$ may contain a point over $0\in\C^2$. We may assume $\Yh$ normal (replacing it by the normalization if necessary) so that it is Cohen-Macaulay by Serre's condition $S_2$, see for instance \cite[II, Theorem 8.22A]{Ha2}. The finiteness and flatness of (A.4.3) then follows from this in the $D_i$ unibranch case (using \cite[II, Theorem 8.21A(c)]{Ha2}). Here $U'$ is a sufficiently small neighborhood in the classical topology and $U$ is a connected component of the inverse image of $U'$. (In the non-unibranch case, etale topology would be needed, see for instance \cite{Ray}.)
\sk
(ii) In the $D_i$ {\it unibranch\1} case, we can also show (A.4.2) by using a remark before (A.4.3) about the coordinate $\widetilde{z}$. Here we may assume $Y$ smooth projective and $D_i\cap D_j=\{y\}$ by replacing $Y$ with a smooth projective compactification of an affine neighborhood of $y\in Y$ and then blowing-up $Y$ at the other intersection points of $D_i,D_j$. Note that the intersection number of the proper transform of $D_i$ and the total transform of $D_j$ coincides with the intersection number $(D_i,D_j)$, since the intersection numbers of the exceptional divisors and the total transform of $D_j$ vanish. So we get (A.4.2) in this case using the remark about $\widetilde{z}$.
\msn
{\bf A.5.~Quotient variety construction.} We can generalize a construction in \cite[Example 6.8]{BR} (due to J.~Koll\'ar and B.~Wang) by using totally ramified cyclic coverings in (A.1) as follows.
\sk
Let $X\to Y$ be the totally ramified cyclic covering of degree $m\ges 2$ associated with $(D,\Lc)$ as in (A.1), where we assume $Y$ connected and smooth as in (A.2). Let $e$ be an integer $\ges 2$ dividing $m$. Let $Z$ be a connected complex algebraic variety having an action of $\muu_e$ such that we have the quotient variety $Z':=Z/\muu_e$ and the canonical morphism $Z\to Z'$ is a finite morphism. (In {\it loc.\,cit.} $Y=\PP^1$, $Z$ is an elliptic curve, $Z\to Z'$ is an isogeny, $m=e=2$, and either $D=2\{0\}$ with $\Lc=\Oc_{\PP^1}(-1)$ or $D=2\{0\}+\{1\}+\{\infty\}$ with $\Lc=\Oc_{\PP^1}(-2)$.)
\sk
There is a short exact sequence of abelian groups
$$1\to\muu_e\to\muu_m\to\muu_{m'}\to 1,$$
where the first morphism is the $m'\1$th power morphism with $m':=m/e$. This gives an action of $\muu_e$ on $X$, and hence the diagonal action of $\muu_e$ on $X\times Z$. We have the quotient variety $(X\times Z)/\muu_e$, since $X$, $Z$, and hence $X\times Z$, are covered by affine open subsets which are stable by the action of $\muu_e$. (Here we use the assertion that a finite morphism is affine.) There is a ramified cyclic covering
$$\XX:=(X\,{\times}\,Z)/\muu_e\,\,\onto\,\,(X\,{\times}\,Z)/(\muu_e\,{\times}\,\muu_e)=(X/\muu_e)\times Z',
\leqno{\rm(A.5.1)}$$
by the diagonal embedding $\muu_e\into\muu_e\,{\times}\,\muu_e$. We can verify that $X/\muu_e$ is the totally ramified cyclic covering of $Y$ of degree $m'=m/e$ associated with $(D,\Lc^{\otimes\,e})$ if $e\ne m$, and $X/\muu_e=Y$ if $e=m$. These imply that $\XX$ is projective if $Y$ is. (Indeed, a finite morphism is projective, see \cite[Propositions 4.4.10 and Corollary 6.1.11]{Gr1}.)
\sk
We have the isomorphisms
$$H^1(\XX,\Q)=H^1(X,\Q)^{\muu_e}\oplus H^1(Z,\Q)^{\muu_e}=H^1(X/\muu_e,\Q)\oplus H^1(Z',\Q).
\leqno{\rm(A.5.2)}$$
The first isomorphism follows from the K\"unneth isomorphism
$$H^1(X\times Z,\Q)=H^1(X,\Q){\otimes}H^0(Z,\Q)\oplus H^0(X,\Q){\otimes}H^1(Z,\Q),
\leqno{\rm(A.5.3)}$$
since the action of $\muu_e$ is trivial on
$$H^0(X,\Q)=H^0(Z,\Q)=\Q.$$
\sk
As a corollary of (A.5.2), we get
$$\h{$W_0H^1(\XX,\Q)=0,\,$ if $\,Z'\,$ is smooth and $\,e=m$,}
\leqno{\rm(A.5.4)}$$
since $X/\muu_e=Y$ is assumed smooth.
\msn
{\bf Remark.} Assume the finite abelian group $\muu_e$ acts on $Z$ freely so that $Z\to Z'$ is an unramified finite covering. For instance, $Z$ is an elliptic curve $E$, or more generally an abelian variety, and the action is given by
$$E\ni Q\to Q+kP\in E\q\h{for}\q k\in\Z/e\Z\,(\cong\muu_e),$$
where $Q$ is a torsion point of $E$ of order $e$.
\sk
In this unramified finite covering case we can get an example where a 1-dimensional stratum $S_i$ in (2.4) is not a rational curve and moreover the monodromy along $S_i$ is nontrivial. Here we have the following isomorphisms compatible with (A.5.2):
$$H^0(Z',R^1(\rho_{\XX})_*\Q_{\XX})=H^1(X,\Q)^{\muu_e}=H^1(X/\muu_e,\Q),
\leqno{\rm(A.5.5)}$$
where $\rho_{\XX}:\XX\to Z'$ is the canonical morphism. For the proof of the first isomorphism, we can verify that the monodromy action on $\bl(R^1(\rho_{\XX})_*\Q_{\XX}\br){}_{z'_0}$ by $\xi\in\pi_1(Z',z'_0)$ is identified with the action of a generator $\gamma\in\muu_e$ on $H^1(X,\Q)$, if $\xi$ is represented by the image of a path from the base point $z_0\in Z$ to $\gamma(z_0)\in Z$. (Note that $\xi$ is unique up to the composition with the image of an element of $\pi_1(Z,z_0)$ in $\pi_1(Z',z'_0)$.) Here $z'_0$ is the base point of $Z'$, which is the image of the base point $z_0\in Z$.
\sk
We have the {\it geometric monodromy} coming from the product structure of $X\times Z$. Indeed, we have the {\it parallel translation} between the fibers
$$\rho_{\XX}^{-1}(z')\cong\rho_{\XX}^{-1}(z'')\,\,\,\,\h{for sufficiently near}\,\,\,z',z''\in Z',
\leqno{\rm(A.5.6)}$$
since we have this for $X\times Z\to Z$ in a compatible way with the action of $\muu_e$. By definition the pull-back of this parallel translation to $Z$ gives the global trivialization of $X\times Z\to Z$. This implies that the geometric monodromy is independent of the choice of a path from $z_0$ to $\gamma(z_0)$ in $Z$ (using the above description of the ambiguity of $\xi\in\pi_1(Z',z'_0)$).
\sk
The second isomorphism of (A.5.5) follows from the isomorphism
$$\bl((\pi_{X'})_*\Q_X\br){}^{\muu_e}=\Q_{X'},
\leqno{\rm(A.5.7)}$$
together with the {\it complete reducibility\1} of finite dimensional complex representations of finite groups (using the scalar extension by $\Q\into\C$), where $\pi_{X'}:X\to X':=X/\muu_e$ is the quotient morphism. This isomorphism is a stalkwise property of the constructible sheaf $(\pi_{X'})_*\Q_X$ with the $\muu_e$-action, since $\pi_{X'}$ is a finite morphism. (This argument holds for any finite group action on $X$ such the quotient exists as an algebraic variety.)


\begin{thebibliography}{RSW}
\bibitem[ADH]{ADH} Arapura, D., Dimca, A.\ and Hain, R., On the fundamental groups of normal varieties, Commun.\ Contemp.\ Math.\ 18 (2016).
\bibitem[BBD]{BBD} Beilinson, A., Bernstein, J.\ and Deligne, P., Faisceaux pervers, Ast\'erisque 100, Soc.\ Math.\ France, Paris, 1982.
\bibitem[BR]{BR} Budur, N.\ and Rubi\'o, M., $L$-infinity pairs and applications to singularities, arXiv:1804.06672.
\bibitem[De1]{De1} Deligne, P., Th\'eorie de Hodge II, Publ.\ Math.\ IHES 40 (1971), 5--57.
\bibitem[De2]{De2} Deligne, P., Le formalisme des cycles \'evanescents, in SGA7 XIII, Lect.\ Notes in Math.\ 340, Springer, Berlin (1973), 82--115.
\bibitem[De3]{De3} Deligne, P., Th\'eorie de Hodge III, Publ.\ Math.\ IHES 44 (1974), 5--77.
\bibitem[Di]{Di} Dimca, A., Sheaves in topology, Universitext, Springer, Berlin, 2004.
\bibitem[Ei]{Ei} Eisenbud, D., Commutative Algebra with a View Toward Algebraic Geometry, Springer, New York, 1994.
\bibitem[Fu]{Fu} Fulton, W., Intersection Theory, Springer, Berlin, 1984.
\bibitem[FL]{FL} Fulton, W.\ and Lazarsfeld, R., Connectivity and its applications in algebraic geometry, in Algebraic Geometry, Lect.\ Notes in Math.\ 862, Springer, Berlin, 1981, pp.~26--92.
\bibitem[GM1]{GM1} Goresky, M.\ and MacPherson, R., Intersection homology theory, Topology 19 (1980), 135--162.
\bibitem[GM2]{GM2} Goresky, M.\ and MacPherson, R., Intersection homology II, Inv.\ Math.\ 71 (1983), 77--129.
\bibitem[Gr1]{Gr1} Grothendieck, A., El\'ements de g\'eom\'etrie alg\'ebrique II, Publ.\ Math.\ IHES 8 (1961).
\bibitem[Gr2]{Gr2} Grothendieck, A., El\'ements de g\'eom\'etrie alg\'ebrique III-1, Publ.\ Math.\ IHES 11 (1961).
\bibitem[Ha1]{Ha1} Hartshorne, R., Ample Subvarieties of Algebraic varieties, Lect.\ Notes in Math.\ 156 Springer, Berlin, 1970.
\bibitem[Ha2]{Ha2} Hartshorne, R., Algebraic Geometry, Springer, Berlin, 1977.
\bibitem[Ko]{Ko} Koll\'ar, J., Shafarevich maps and automorphic forms, Princeton Univ.\ Press, Princeton, NJ, 1995.
\bibitem[Ma]{Ma} Mather, J., Notes on topological stability (Harvard, 1970), Bull.\ Amer.\ Math.\ Soc.\ (N.S.) 49 (2012), 475--506.
\bibitem[McC]{McC} McCrory, C., Poincar\'e duality in spaces with singularities, Ph.D.\ Thesis, Brandeis University (1972).
\bibitem[Mi]{Mi} Milnor, J., Singular Points of Complex Hypersurfaces, Princeton Univ.\ Press, Princeton, NJ, 1968.
\bibitem[Mu]{Mu} Mumford, D., The red book of varieties and schemes, Lect.\ Notes in Math. 1358, Springer, Berlin, 1999.
\bibitem[Na]{Na} Navarro Aznar, V., Th\'eor\`emes d'annulation, Lect.\ Notes in Math.\ 1335, Springer, Berlin, 1988, pp.~133--160.
\bibitem[PS]{PS} Peters, C.A.M.\ and Steenbrink, J.H.M., Mixed Hodge structures, Springer, Berlin, 2008.
\bibitem[Ram]{Ram} Ramanujam, C.P., Remarks on the Kodaira vanishing theorem, J.\ Indian Math.\ Soc.\ (N.S.) 36 (1972), 41--51.
\bibitem[Ray]{Ray} Raynaud, M., Anneaux locaux hens\'eliens, Lect.\ Notes in Math.\ 169, Springer, Berlin, 1970.
\bibitem[RSW]{RSW} Reichelt, T., Saito, M.\ and Walther, U., Dependence of Lyubeznik numbers of cones of projective schemes on projective embeddings, arXiv:1803.07448.
\bibitem[Sa1]{mhp} Saito, M., Modules de Hodge polarisables, Publ. RIMS, Kyoto Univ. 24 (1988), 849--995.
\bibitem[Sa2]{mhm} Saito, M., Mixed Hodge modules, Publ. RIMS, Kyoto Univ.\ 26 (1990), 221--333.
\bibitem[Sa3]{mhc} Saito, M., Mixed Hodge complexes on algebraic varieties, Math.\ Ann.\ 316 (2000), 283--331.
\bibitem[Sa4]{wh} Saito, M., Bernstein-Sato polynomials for projective hypersurfaces with weighted homogeneous isolated singularities, arXiv:1609.04801.
\bibitem[Se]{Se} Serre, J.-P., Algebre locale, multiplicit\'es, Lect.\ Notes in Math.\ 11, Springer, Berlin, 1975.
\bibitem[St1]{St1} Steenbrink, J.H.M., Intersection form for quasi-homogeneous singularities, Compos.\ Math.\ 34 (1977), 211--223.
\bibitem[St2]{St2} Steenbrink, J.H.M., Mixed Hodge structure on the vanishing cohomology, in Real and complex singularities, Sijthoff and Noordhoff, Alphen aan den Rijn, 1977, pp. 525--563.
\bibitem[SS]{SS} Steenbrink, J.H.H.\ and Stevens, J., Topological invariance of weight filtration, Indag.\ Math.\ 46 (1984), 63--76.
\bibitem[Ve]{Ve} Verdier, J.-L., Dualit\'e dans la cohomologie des espaces localement compacts, S\'eminaire Bourbaki, Vol. 9, Exp.\ No.\ 300, Soc.\ Math.\ France, Paris, 1995, pp.~337--349.
\bibitem[Vi]{Vi} Viehweg, E., Vanishing theorems, J.\ reine angew.\ Math.\ 335 (1982), 1--8.
\bibitem[We]{We} Weber, A., Pure homology of algebraic varieties, Topology 43 (2004), 635--644.
\end{thebibliography}
\end{document}